\documentclass{amsart}
\usepackage{graphicx}

\newtheorem{thm}{Theorem}[section]

\theoremstyle{definition}

\theoremstyle{remark}

\numberwithin{equation}{section}

\begin{document}

\title[Continuity of Large Closed Queueing Networks with Bottlenecks]{Continuity of Large Closed Queueing Networks with Bottlenecks}%
\author{Vyacheslav M. Abramov}%
\address{School of Mathematical Sciences, Monash University,
Clayton Campus, Wellington road, VIC-3800, Australia}%
\email{Vyacheslav.Abramov@sci.monash.edu.au}%

\thanks{The author acknowledges with thanks the support of the Australian
Research Council}%

\begin{abstract}
This paper studies a closed queueing network containing a hub (a
state dependent queueing system with service depending on the number
of units residing here) and $k$ satellite stations, which are
$GI/M/1$ queueing systems. The number of units in the system, $N$,
is assumed to be large. After service completion in the hub, a unit
visits a satellite station $j$, $1\leq j\leq k$, with probability
$p_j$, and, after the service completion there, returns to the hub.
The parameters of service times in the satellite stations and in the
hub are proportional to $\frac{1}{N}$. One of the satellite stations
is assumed to be a bottleneck station, while others are
non-bottleneck. The paper establishes the continuity of the
queue-length processes in non-bottleneck satellite stations of the
network when the service times in the hub are close in certain sense
(exactly defined in the paper) to the exponential distribution.
\end{abstract}
\maketitle
\newpage

\section{Introduction}\label{Introduction}

In the present paper we study the continuity of large closed
queueing networks with bottlenecks. The continuity of queueing
systems is a known area in queueing theory, and the goal of this
theory is to answer the question: \textit{how small change of input
characteristics of the system affects its output characteristics}.
This question is of especial significance for queueing networks that
models telecommunication networks.

In the present paper the continuity results are established for
closed queueing networks describing client/server telecommunication
networks. The queueing networks considered in the paper are of the
following configuration. There is a large server (hub), which is a
specific state dependent queueing system (the details are given
later) and $k$ satellite stations, which are single server
$GI/M/1/\infty$ queueing systems. Being served in the hub, a unit is
addressed to the $j$th satellite station, $1\leq j\leq k$, with
probability $p_j>0$ $\left(\sum\limits_{j=1}^{k}p_j=1\right)$, and
then, being served there, it returns to the hub.

There are many papers in the literature that study similar type of
queueing systems. In most of them the hub is a multiserver (or
infinite server) queueing system. The earliest consideration seems
to go back to a paper of Whitt \cite{Whitt 1984}. A Markovian
queueing network with one bottleneck satellite station has been
studied by Kogan and Liptser \cite{Kogan and Liptser 1993}.
Extensions of the results \cite{Kogan and Liptser 1993} for
different non-Markovian models have been made in a series of papers
and a book by the author (see \cite{Abramov 2000}, \cite{Abramov
2001}, \cite{Abramov 2004}, \cite{Abramov 2008} and \cite{Abramov
2009}). Non-Markovian models of networks with a hub and one
satellite station have also been considered by Krichagina, Liptser
and Puhalskii \cite{Krichagina Liptser and Puhalsky 1989} and
Krichagina and Puhalskii \cite{Krichagina and Puhalsky 1998}. For
other studies of queueing networks with bottleneck see Berger,
Bregman and Kogan \cite{Berger Bregman and Kogan 1999}, Brown and
Pollett \cite{Brown and Pollett 1982} and Pollett \cite{Pollett
2000} and other papers.

Aforementioned papers \cite{Abramov 2000}, \cite{Abramov 2001},
\cite{Abramov 2004}, \cite{Abramov 2008} and \cite{Kogan and Liptser
1993}, study non-stationary queue-length distributions in
non-bottleneck satellite stations. Specifically, in \cite{Abramov
2000}, \cite{Abramov 2004} and \cite{Kogan and Liptser 1993} these
distributions have been studied in the presence of one bottleneck
satellite station, and in \cite{Abramov 2001} the only cases of
bottleneck leading to diffusion approximations have been
investigated. One of the main assumption leading to substantial
technical simplification in \cite{Abramov 2000}, \cite{Abramov
2001}, \cite{Abramov 2004} and \cite{Kogan and Liptser 1993} was
that at the initial time moment all of $N$ units are located in the
hub, where the time instant $t=0$ is the moment of service starts
for all of them. In the case of exponentially distributed service
times in the hub, because of the property of the lack of memory of
an exponential distribution, the above assumption is not the loss of
generality. The initial time moment $t=0$ is not necessarily the
moment of service starts for the units presented at the hub, and the
number of these units need not be $N$ exactly; it can be a large
random number $\mathcal{N}_0$ which is asymptotically equivalent to
$N$, i.e. $\frac{\mathcal{N}_0}{N}$ converges in distribution to 1
as $N$ increases unboundedly.

However, in the case where the service times in the hub are
generally distributed, the above assumption is an essential
restriction, because its change can lead to a serious change of the
limiting non-stationary queue-length distributions. The models
considered in \cite{Abramov 2001} (see also \cite{Abramov 2009},
Chap. 2) are flexible to this change, since they model hub as a
state dependent single-server queueing system. The behavior of the
queue-length processes in that system has been studied analytically,
and the analysis of that system in \cite{Abramov 2001} and
\cite{Abramov 2009} shows that it is stable with respect to a slight
change of the initial condition. The assumption on state dependence
in \cite{Abramov 2001} and \cite{Abramov 2009}, Chap. 2,  was quite
general, and real bottleneck analysis for that network, except the
cases where diffusion approximations are available, could not be
provided. So, the most interesting cases of the problem of the
bottleneck analysis in \cite{Abramov 2001} or \cite{Abramov 2009},
Chap. 2 remained unsolved.


In the present paper we study a more particular model than that in
\cite{Abramov 2001}, i.e. we study the model with the same
state-dependent service mechanism in the hub as in \cite{Abramov
2001}, but specific additional assumptions on probability
distribution functions of service times in the hub that will be
discussed later in the paper. One of these assumptions, saying that
service times in the hub are ``close" (in different senses given in
the paper) to the exponential distribution, enables us to easily
extend the results to various models, including the realistic cases
where the hub is multiserver (or infinite-server) queueing system.

%

Recall the assumptions that made in \cite{Abramov 2001} in the
description of the service mechanism in the hub. In \cite{Abramov
2001} it was assumed as follows.

\smallskip
\noindent ({\bf i}) if immediately before a service start of a unit
the queue-length in the hub is equal to $K\leq N$, then the
probability distribution function is $F_K(x)=G_K(Kx)$. It was
assumed that $\frac{1}{\lambda_K}=\int_0^\infty x\mathrm{d}G_K(x)$.

\smallskip
Hence, for the expectation of a service time we have
$$
\int_0^\infty x\mathrm{d}F_K(x)=\int_0^\infty
x\mathrm{d}G_K(Kx)=\frac{1}{K}\int_0^\infty
x\mathrm{d}G_K(x):=\frac{1}{K\lambda_K}.
$$

\noindent ({\bf ii}) As $N\to\infty$, the sequence of probability
distribution functions $G_N(x)$ was assumed to converge weakly to
$G(x)$ with $G(0+)=0$, and $\frac{1}{\lambda}=\int_0^\infty
x\mathrm{d}G(x)$. In addition, it was assumed in \cite{Abramov 2001}
that $G_K(Kx)\leq G_{K+1}(Kx+x)$ for all $x\geq 0$.

\smallskip
In the present paper in addition to ({\bf i}) we assume that
\begin{equation}\label{0.1add}
G_1(x)=G_2(x)=\ldots=G_N(x):=G(x),
\end{equation}
and there exists the second moment
$$
\int_0^\infty x^2\mathrm{d}G(x)=r<\infty
$$
satisfying the condition:
\begin{equation}\label{f1}
r>\frac{2}{\lambda^{2}}.
\end{equation}
Condition \eqref{f1} will be used later in the proof of Theorem
\ref{thm1}.

In this case, for the random variable $\zeta_K$ having the
probability distribution function $F_K(x)=G(Kx)$ we have
\begin{equation}\label{0.1add.2}
\mathrm{E}\zeta_K=\int_0^\infty x \mathrm{d}F_K(x)=\int_0^\infty x
\mathrm{d}G(Kx)=\frac{1}{K\lambda}.
\end{equation}
Under the above addition assumptions, earlier assumption ({\bf ii})
of \cite{Abramov 2001} becomes irrelevant here.

Following \eqref{0.1add} and \eqref{0.1add.2} one can ``guess" that
departures of units from the hub when $K$ increases indefinitely
should be asymptotically insensitive to the type of the distribution
$G(x)$, i.e. the number of departures per unit of time should be the
same as in the case when $G(x)$ would be the exponential
distribution. This guess however cannot be proved by the means of
the earlier papers \cite{Kogan and Liptser 1993} and \cite{Abramov
2000} and should involve more delicate methods (which are provided
in the present paper).
\smallskip

In addition to the assumptions made, we consider two classes of
probability distribution function $G(x):=\mathrm{P}\{\zeta\leq x\}$.
The first one satisfies the condition
\begin{equation}\label{0.2}
\sup_{x,y\geq 0}|G(x)-G_y(x)|<\epsilon,
\end{equation}
where $\epsilon$ is a small positive number, and
$G_y(x)=\mathrm{P}\{\zeta\leq x+y|\zeta>y\}$.

For the second class of distributions along with (\ref{0.2}) it is
assumed that $G(x)$ belongs to one of the classes NBU (New Better
than Used) or NWU (New Worse that Used). Recall that a probability
distribution function $\Xi(x)$ of a nonnegative random variable is
said to belong to the class NBU if for all $x\geq0$ and $y\geq0$ we
have $\overline{\Xi}(x+y)\leq\overline{\Xi}(x)\overline{\Xi}(y)$,
where $\overline{\Xi}(x)=1-{\Xi}(x)$. If the opposite inequality
holds, i.e.
$\overline{\Xi}(x+y)\geq\overline{\Xi}(x)\overline{\Xi}(y)$, then
$\Xi(x)$ is said to belong to the class NWU.

The assumption given by \eqref{0.2} involves many nice properties to
the probability distribution function $G(x)$. We review these
properties now, but the details will be given then.

First, according to \eqref{0.2} the probability distribution $G(x)$
is close to the exponential distribution in the uniform (Kolmogorov)
metric in the sense that if $\epsilon$ is small, then the distance
in Kolmogorov's metric between $G(x)$ and the exponential
distribution with mean $\frac{1}{\lambda}$ is small as well (see
later relations \eqref{0.5} and \eqref{0.7}).

Second, according to assumption ({\bf i}) the probability
distribution of a current service time in the hub is not changed
when a new arrival occurs from one of satellite stations to the hub
during that service time. The change of service time distribution in
the hub can be at the moments of service start only. Additional
assumption \eqref{0.2} enables us to consider more general schemes
in which the change of probability distribution function of a
service time in the hub can be taken into account also at the
moments of arrival of customers to the hub from satellite stations.
This means that the present model can be essentially extended, and
the extension can include the models where the hub is a multiserver
or infinite-server queueing system.

The details for the explanation of the second property are given in
\eqref{nice2.1} and \eqref{nice2.2} below. Specifically, under
assumption \eqref{0.2} we have:

\begin{eqnarray}
&&\sup_{x,y_1\geq 0}|G_{y_1}(x)-G_{y_2}(x)|\label{nice2.1}\\
&&\leq\sup_{x,y_1\geq 0}|G(x)-G_{y_1}(x)|+\sup_{x\geq
0}|G(x)-G_{y_2}(x)|\nonumber\\ &&\leq\sup_{x,y_1\geq
0}|G(x)-G_{y_1}(x)|+\sup_{x,y_2\geq 0}|G(x)-G_{y_2}(x)|\nonumber\\
&&<2\epsilon.\nonumber
\end{eqnarray}
Therefore
\begin{equation}\label{0.4}
\sup_{x,y_1,y_2\geq 0}|G_{y_1}(x)-G_{y_2}(x)|<2\epsilon
\end{equation}
as well. As well, according to the characterization theorem of
Azlarov and Volodin \cite{Azlarov and Volodin 1987}, \cite{Abramov
2008b} we have:
\begin{equation}\label{0.5}
\sup_{x\geq 0}|G(x)-(1-\mathrm{e}^{-\lambda x})|<2\epsilon.
\end{equation}

If, in addition, the probability distribution function $G(x)$
belongs either to the class NBU or to the class NWU, then, instead
of (\ref{0.4}), for any $y_2\geq0$ we, respectively, have:
\begin{eqnarray}
\sup_{x,y_1\geq 0}|G_{y_1}(x)-G_{y_2}(x)|
&&\leq\sup_{x,y_1\geq 0}\max_{i=1,2}|G(x)-G_{y_i}(x)|\label{nice2.2}\\
&&\leq\sup_{x,y_1\geq
0}|G(x)-G_{y_1}(x)|\nonumber\\
&&<\epsilon.\nonumber
\end{eqnarray}

Hence,
\begin{equation}\label{0.6}
\sup_{x,y_1,y_2\geq 0}|G_{y_1}(x)-G_{y_2}(x)|<\epsilon.
\end{equation}
As well, instead of (\ref{0.5}) (see \cite{Abramov 2008b}) we have
\begin{equation}\label{0.7}
\sup_{x\geq 0}|G(x)-(1-\mathrm{e}^{-\lambda x})|<\epsilon.
\end{equation}

As in \cite{Abramov 2001}, it is assumed in the paper that the
service times in the satellite station $j$, $j=1,2,\ldots,k$, are
exponentially distributed with parameter $N\mu_j$, and
$\mu_j>\lambda p_j$ for $j=1,2,\ldots, k-1$ and $\mu_k<\lambda p_k$.
This means that the $k$th satellite station is assumed to be a
bottleneck station, while the first $k-1$ satellite stations are
non-bottleneck. The only this case is studied in the paper.

The problem of continuity of queues and networks is an old problem.
It goes back to the papers of Kennedy \cite{Kennedy 1972} and Whitt
\cite{Whitt 1974}, and nowadays there are many papers in this area.
To mention only a few of them we refer Kalashnikov
\cite{Kalashnikov}, Kalashnikov and Rachev \cite{Kalashnikov and
Rachev}, Rachev \cite{Rachev}, Zolotarev \cite{Z1}, \cite{Z2},
Gordienko and Ruiz de Ch\'avez \cite{GordR-Ch1}, \cite{GordR-Ch2}.

 The continuity results of the present paper are based on the
estimates for the least positive root of the functional equation
\begin{equation}\label{TakEq}
z=\widehat{G}(\mu-\mu z),
\end{equation}
where $\widehat{G}(s)$, $s\geq0$, is the Laplace-Stieltjes transform
of the probability distribution function $G(x)$ of a positive random
variable having first two moments, and satisfying the condition
$\mu\int_0^\infty x\mathrm{d}G(x)>1$. It is well-known (see Tak\'acs
\cite{Takacs 1962} or Gnedenko and Kovalenko \cite{Gnedenko and
Kovalenko 1989}) that the least positive root of the aforementioned
functional equation \eqref{TakEq} belongs to the interval (0,1).
Within the interval (0,1) Rolski \cite{Rolski 1972} established the
possible lower and upper bounds for the least positive root of
equation \eqref{TakEq}, when the probability distribution function
$G(x)$ has two moments. Let $\mathcal{G}(\frak{g}_1, \frak{g}_2)$ be
the class of probability distribution functions of positive random
variables having the first two moments $\frak{g}_1$ and
$\frak{g}_2$. Let $G^{(1)}(x)$ and $G^{(2)}(x)$ be two probability
distribution functions from this class satisfying the additional
condition $\sup_{x\geq0}|G^{(1)}(x)-G^{(2)}(x)|<\epsilon$. Let
$\widehat G^{(1)}(s)$ and $\widehat G^{(2)}(s)$, $s\geq0$, denote
the Laplace-Stieltjes transforms of $G^{(1)}(x)$ and $G^{(2)}(x)$
respectively, and let $\gamma^{(1)}$ and $\gamma^{(2)}$ denote the
corresponding least roots of the functional equations
$z=\widehat{G}^{(1)}(\mu-\mu z)$ and $z=\widehat{G}^{(2)}(\mu-\mu
z)$. The lower and upper bounds for $|\gamma^{(1)}-\gamma^{(2)}|$
have been established in \cite{Abramov cont}. In the present paper
we aimed to use these bounds to establish the desired continuity
results for the nonstationary queue-length distributions in the
non-bottleneck stations of our network.

The rest of the paper is organized as follows. In Section
\ref{OneSatellite} the known results from \cite{Kogan and Liptser
1993} for Markovian models in which the hub is infinite-server
queueing system and there is only one satellite station, which is a
bottleneck station, are recalled. In same Section
\ref{OneSatellite}, we prove that the same result holds true under
the main assumption of this paper where the service times in the hub
are state dependent and satisfy \eqref{0.1add}. In Section
\ref{General}, we extend our result to the case of our model where
there are \textit{several} satellite stations and one of them is
bottleneck. In Section \ref{Continuity}, the main results of this
paper related to the continuity of non-stationary queue-length
distributions in non-bottleneck satellite stations are obtained.
Some technical auxiliary results are recalled in the Appendix.

\section{A large closed network containing only one satellite
station, which is bottleneck}\label{OneSatellite} In this section we
consider examples of simplest queueing networks in order to
demonstrate that in the case of the model considered in the paper,
the asymptotic behaviour of the queue-length process in the hub is
the same as that in the case of the exponentially distributed
service times in the hub, which was originally studied in
\cite{Kogan and Liptser 1993}.

\subsection{The case $G(x)=1-\mathrm{e}^{-\lambda x}$}
In the case where $G(x)=1-\mathrm{e}^{-\lambda x}$ we deal with the
model considered in the paper of Kogan and Liptser \cite{Kogan and
Liptser 1993}.

Consider the Markovian network (see Kogan and Liptser \cite{Kogan
and Liptser 1993}) containing only one satellite station, which is
assumed to be bottleneck. Let $A(t)$ and $D(t)$ denote arrival
process to this satellite station and, respectively, departure
process from this satellite station. Denoting the queue-length
process in this satellite station by $Q(t)$ one can write:
\begin{equation}\label{1.1}
Q(t)=A(t)-D(t).
\end{equation}
The departure process $D(t)$ is defined via the Poisson process as
follows. Service times in the satellite station are assumed to be
exponentially distributed with parameter $\mu N$. Therefore,
denoting the Poisson process with parameter $\mu N$ by $S(t)$, for
the departure process $D(t)$ we have the equation:
\begin{equation}\label{1.2}
D(t)=\int_0^t\mathbf{1}\{Q(s-)>0\}\mathrm{d}S(s),
\end{equation}
where $\mathbf{1}\{\mathcal{A}\}$ denotes the indicator of the
event $\mathcal{A}$.

To define the arrival process, consider the collection of
independent Poisson processes $\pi_i(t)$, $i=1,2,\ldots,N$. Then,
\begin{equation}\label{1.3}
A(t)=\int_0^t\sum_{i=1}^{N}\mathbf{1}\{N-Q(s-)\geq
i\}\mathrm{d}\pi_i(s).
\end{equation}

For the normalized queue-length process $q_N(t)=\frac{1}{N}Q(t)$,
from (\ref{1.1}) we have:
\begin{equation}\label{1.4}
q_N(t)=\frac{1}{N}A(t)-\frac{1}{N}D(t).
\end{equation}
As $N\to\infty$, both $\frac{1}{N}A(t)$ and $\frac{1}{N}D(t)$
converge a.s. to the limits, and these limits are the same as the
correspondent limits of $\frac{1}{N}\widehat{A}(t)$ and
$\frac{1}{N}\widehat{D}(t)$ (see \cite{Kogan and Liptser 1993} or
\cite{Abramov 2000} for further details), where $\widehat{A}(t)$
and $\widehat{D}(t)$ are the compensators in the Doob-Meyer
semimartingale decompositions $A(t)=M_A(t)+\widehat{A}(t)$ and,
correspondingly, $D(t)=M_D(t)+\widehat{D}(t)$ ($M_A(t)$ and
$M_D(t)$ denote local square integrable martingales corresponding
the processes $A(t)$ and $D(t)$.) Therefore,
\begin{equation}\label{1.5}
q_N(t)=\frac{1}{N}\widehat{A}(t)-\frac{1}{N}\widehat{D}(t).
\end{equation}
The compensators $\widehat{A}(t)$ and $\widehat{D}(t)$ have the
representations (for details see \cite{Kogan and Liptser 1993}):
\begin{equation}\label{1.6}
\widehat{A}(t)=\lambda\int_0^t[N-Q(s)]\mathrm{d}s,
\end{equation}
and
\begin{equation}\label{1.7}
\widehat{D}(t)=\mu Nt-\mu
N\int_0^t\mathbf{1}\{Q(s)=0\}\mathrm{d}s,
\end{equation}

Therefore, from (\ref{1.5}), (\ref{1.6}) and (\ref{1.7})
\begin{eqnarray}\label{1.8}
q(t):&=&\mathrm{P}^\_\lim_{N\to\infty}q_N(t)\\
&=&\lambda\int_0^t[1-q(s)]\mathrm{d}s -\mu
t\nonumber\\
&&+\mu\int_0^t\mathbf{1}\{Q(s)=0\}\mathrm{d}s.\nonumber
\end{eqnarray}
In the case $\mu<\lambda$, using the Skorokhod reflection
principle and techniques of the theory of martingales one can show
that $\mathbf{1}\{Q(s)=0\}$ vanishes in probability for all $s>0$
as $N$ increases to infinity (for details of this see \cite{Kogan
and Liptser 1993} or \cite{Abramov 2000}), and from (\ref{1.8}) we
arrive at
\begin{equation}\label{1.9}
q(t)=\lambda\int_0^t[1-q(s)]\mathrm{d}s-\mu t.
\end{equation}
Taking into account that $q(0)=1$, from (\ref{1.9}) we obtain:
\begin{equation}\label{1.10}
q(t)=\left(1-\frac{\mu}{\lambda}\right)(1-\mathrm{e}^{-\lambda
t}).
\end{equation}

\subsection{The case of arbitrary distribution $G(x)$}
We will show below that the equation (\ref{1.10}) holds also in the
case, when service times in the hub are generally distributed and
state dependent (for further details see Section \ref{Introduction})
but under the assumption that there is only one satellite station,
which is assumed to be bottleneck.

Let $\delta$ be a small enough positive fixed real number so the
inequalities appearing later in the form $\mu<\lambda(1-c\delta)$
for specified constants $c$ are assumed to be satisfied; the large
parameter $N$ is assumed to increase to infinity. Let $Q_N(t)$
denote the queue-length in the satellite station (which is a
bottleneck station), and let $\overline{Q}_N(t)$ denote the
queue-length in the hub ($\overline{Q}_N(t)=N-Q_N(t)$). In a small
interval $[0,\delta)$ let us denote
$$\overline{Q}_N^+[0,\delta)=\sup_{0\leq
t<\delta}\overline{Q}_N(t)$$ and
$$\overline{Q}_N^-[0,\delta)=\inf_{0\leq
t<\delta}\overline{Q}_N(t).$$ Clearly, that
$\overline{Q}_N^+[0,\delta)=N$, and, as $N\to\infty$, the a.s. lower
and upper limits of $\overline{q}_N(t)=\frac{1}{N}\overline{Q}_N(t)$
as $N\to\infty$, must be between 1 and $1-(\lambda-\mu)\delta$ for
all $t\in[0,\delta)$.

The last result is explained as follows. Let $a_N[0,\delta)$ be the
maximum arrival rate (that is maximum instantaneous arrival rate) in
the interval $[0,\delta)$. Since $\overline{Q}_N^+[0,\delta)=N$,
this instantaneous rate $a_N[0,\delta)=\lambda N$. This means that,
as $N\to\infty$, then with probability approaching 1 the number of
arrivals during the interval $[0, \delta)$ is not greater than
$\lambda N \delta$, i.e.
\begin{equation}\label{a1}
\mathrm{Pr}\left\{\limsup_{N\to\infty}\frac{\mathcal{A}_N[0,\delta)}{N}\leq\lambda\delta\right\}=1,
\end{equation}
where $\mathcal{A}_N[0,\delta)$ denotes the number of unit arrivals
from the hub to satellite station. We use $\limsup_{N\to\infty}$
rather than $\lim_{N\to\infty}$ because interarrival times are not
identically distributed, so existence of the limit is not evident.

Let $\mathcal{S}_N[0,\delta)$ denote the number of service
completions in the satellite station during the interval
$[0,\delta)$. Similarly to \eqref{a1}, as $N\to\infty$, we have
\begin{equation}\label{a2}
\mathrm{Pr}\left\{\lim_{N\to\infty}\frac{\mathcal{S}_N[0,\delta)}{N}=\mu\delta\right\}=1,
\end{equation}
which means that, as $N\to\infty$, then with the probability
approaching 1 the number of service completions during [0, $\delta$)
becomes close to $\mu N\delta$.

The fact given in \eqref{a2} can be easily proved as follows. Take
$\epsilon=\frac{\delta}{M}$, where $M$ is a sufficiently large
number. As $\mu<\lambda$, the probability that a busy period of the
$GI/M/1/\infty$ queueing system is finite is $\frac{\mu}{\lambda}$
(i.e. it is less than 1). Therefore the value $N$ can be chosen such
that during the time interval [0, $\epsilon$) there is only a finite
number of busy periods with probability 1, and during the time
interval [$\epsilon$, $\delta$) there is only an unfinished busy
period, and units are served continuously without delay with the
rate $\mu N$ per time unit. Therefore, as $N\to\infty$, with
probability approaching 1 the number of service completions during
the interval [0, $\delta$), $\mathcal{S}_N[0,\delta)$, is at least
$\mu N (\delta-\epsilon)=\mu N\frac{M-1}{M}\delta$. As $N\to\infty$,
the sequence of values $M=M(N)$ can be taken increasing to infinity
and associated sequence $\epsilon=\epsilon(N)$ vanishing. Therefore,
as $N\to\infty$, the value $\mathcal{S}_N[0,\delta)$ becomes
asymptotically equivalent to $\mu N\delta$. So, we arrive at
\eqref{a2}. It also follows from \eqref{a1} and \eqref{a2} that
\begin{equation}\label{a1.1}
\mathrm{Pr}\left\{\liminf_{N\to\infty}
\frac{\mathcal{A}_N[0,\delta)}{N}\geq\lambda[1-(\lambda-\mu)\delta]\delta\right\}=1.
\end{equation}

From \eqref{a1} and \eqref{a2} we obtain
\begin{equation}\label{a3}
\begin{aligned}
&\mathrm{Pr}\left\{1-\limsup_{N\to\infty}\left(\frac{\mathcal{A}_N[0,\delta)}{N}
-\frac{\mathcal{S}_N[0,\delta)}{N}\right)\leq\liminf_{N\to\infty}\overline{q}_{N}(t)\right.\\
& \ \ \ \left.\leq\limsup_{N\to\infty}\overline{q}_{N}(t)\leq1\right\}\\
&=\mathrm{Pr}\left\{1-(\lambda-\mu)\delta\leq\liminf_{N\to\infty}\overline{q}_{N}(t)
\leq\limsup_{N\to\infty}\overline{q}_{N}(t)\leq1\right\}=1.
\end{aligned}
\end{equation}
Hence, according to \eqref{a3}
\begin{equation}\label{a4}
1-(\lambda-\mu)\delta \ {\buildrel{a.s}\over \leq}\overline{q}(t) \
 {\buildrel{a.s}\over \leq} \ 1
\end{equation}
for all $t\in[0,\delta)$, where
$\overline{q}(t)=\frac{1}{2}\left[\liminf_{N\to\infty}\overline{q}_{N}(t)
+\limsup_{N\to\infty}\overline{q}_{N}(t)\right]$. Note, that as
$\delta$ vanishes, $\liminf_{N\to\infty}\overline{q}_{N}(t)$ and
$\limsup_{N\to\infty}\overline{q}_{N}(t)$ approach one another (see
rel. \eqref{a1} and \eqref{a1.1}), and from \eqref{a4} we also have
\begin{equation}\label{a4.1}
\lim_{\delta\downarrow0}\frac{1-\overline{q}(\delta)}{\delta}=\lambda-\mu.
\end{equation}
So, the value $\delta$ can be chosen such small that
$\overline{q}(\delta)$ is asymptotically close to
$1-(\lambda-\mu)\delta$. The existence of the similar limit is also
obtained later in \eqref{a4.2}, \eqref{a4.3} and \eqref{1.11}.


Let us consider the second interval [$\delta$, $2\delta$). As $N$ is
large enough, for small $\delta$ the number of units in the hub at
time $\delta$, $\overline{Q}_N(\delta)$, behaves as
$N[1-(\lambda-\mu)\delta]$, i.e we have approximately this number of
busy servers in the hub. Therefore, at that time moment $\delta$ the
instantaneous rate of arrival of units in the satellite station is
asymptotically equivalent to $\lambda N[1-(\lambda-\mu)\delta]$ and
the service rate $\mu N$. Assuming the inequality
$\mu<\lambda[1-(\lambda-\mu)\delta]$ satisfied, one can use the same
arguments as above. First, $a_N[\delta,2\delta)$ is asymptotically
equivalent to $\lambda N[1-(\lambda-\mu)\delta]$. As $N\to\infty$,
\begin{equation}\label{a5}
\mathrm{Pr}\left\{\limsup_{N\to\infty}\frac{\mathcal{A}_N[\delta,2\delta)}{N}\leq
\lambda\delta[1-(\lambda-\mu)\delta]\right\}=1,
\end{equation}
and
\begin{equation}\label{a6}
\mathrm{Pr}\left\{\lim_{N\to\infty}\frac{\mathcal{S}_N[\delta,2\delta)}{N}=\mu\delta\right\}=1,
\end{equation}
and at the moment $2\delta$ the queue-length in the satellite
station, $Q_N(2\delta)$, is asymptotically evaluated as
$\frac{\lambda-\mu}{\lambda}[2(\lambda\delta)-(\lambda\delta)^{2}]N$.
It also follows from \eqref{a5} and \eqref{a6} that
\begin{equation}\label{a5.1}
\mathrm{Pr}\left\{\liminf_{N\to\infty}\frac{\mathcal{A}_N[\delta,2\delta)}{N}\geq
\lambda\delta\left[1-\frac{\lambda-\mu}{\lambda}[2(\lambda\delta)-(\lambda\delta)^{2}]\right]\right\}=1,
\end{equation}

From \eqref{a5} and \eqref{a6} we obtain:
\begin{equation}\label{a7}
\begin{aligned}
&\mathrm{Pr}\left\{1-\limsup_{N\to\infty}\left(\frac{\mathcal{A}_N[0,2\delta)}{N}
-\frac{\mathcal{S}_N[0,2\delta)}{N}\right)\leq\liminf_{N\to\infty}\overline{q}_{N}(t)\right.\\
&\ \ \
\left.\leq\limsup_{N\to\infty}\overline{q}_{N}(t)\leq1-\limsup_{N\to\infty}\left(\frac{\mathcal{A}_N[0,\delta)}{N}
-\frac{\mathcal{S}_N[0,\delta)}{N}\right)\right\}\\
&=\mathrm{Pr}\left\{1-\frac{\lambda-\mu}{\lambda}[2(\lambda\delta)-(\lambda\delta)^{2}]
\leq\liminf_{N\to\infty}\overline{q}_{N}(t)\right.\\
&\ \ \
\left.\leq\limsup_{N\to\infty}\overline{q}_{N}(t)\leq1-(\lambda-\mu)\delta\right\}=1.
\end{aligned}
\end{equation}
Hence, according to \eqref{a7}
\begin{equation}\label{a8}
1-\frac{\lambda-\mu}{\lambda}[2(\lambda\delta)-(\lambda\delta)^{2}]
\ {\buildrel{a.s}\over \leq}\overline{q}(t) \
 {\buildrel{a.s}\over \leq} \ 1-(\lambda-\mu)\delta
\end{equation}
for all $t\in[\delta,2\delta)$. As $\delta$ vanishes, similarly to
\eqref{a4.1} we have
\begin{equation}\label{a4.2}
\lim_{\delta\downarrow0}\frac{\overline{q}(2\delta)-2\overline{q}(\delta)+1}{\delta^2}=\lambda(\lambda-\mu).
\end{equation}

Then, considering the third interval [$2\delta$, $3\delta$), at the
end of this interval the queue-length in the satellite station,
$Q_N(3\delta)$, is evaluated as
$\frac{\lambda-\mu}{\lambda}[3(\lambda\delta)-3(\lambda\delta)^{2}+(\lambda\delta)^3]N$.
In this case, similarly to \eqref{a8} we arrive at the inequality
\begin{equation}\label{a9}
1-\frac{\lambda-\mu}{\lambda}[3(\lambda\delta)-3(\lambda\delta)^{2}+(\lambda\delta)^3]
\ {\buildrel{a.s}\over \leq}\overline{q}(t) \
 {\buildrel{a.s}\over \leq} \ 1-\frac{\lambda-\mu}{\lambda}[2(\lambda\delta)-(\lambda\delta)^{2}]
\end{equation}
for all $t\in[2\delta,3\delta)$. As $\delta$ vanishes, similarly to
\eqref{a4.1} and \eqref{a4.2},
\begin{equation}\label{a4.3}
\lim_{\delta\downarrow0}\frac{\overline{q}(3\delta)-3\overline{q}(2\delta)+3\overline{q}(\delta)-1}{\delta^{3}}
=\lambda^{2}(\lambda-\mu).
\end{equation}

Thus, considering the sequence of intervals, at the end of the $j$th
interval, for the queue-length in the satellite station,
$Q_N(j\delta)$, we obtain the expansion:
\begin{eqnarray}
&&\frac{\lambda-\mu}{\lambda}\left[j(\delta\lambda)-\left(\begin{array}{c}j\\
2\end{array}\right)(\delta\lambda)^2+\ldots\right.\nonumber\\
&&+(-1)^i\left(\begin{array}{c}j\\
i\end{array}\right)(\delta\lambda)^{i+1}
+\ldots\nonumber\\
&&+(-1)^{j}\left.(\delta\lambda)^{j}\right]\nonumber.
\end{eqnarray}
From this expansion we have
\begin{equation}\label{a11}
\begin{aligned}
&1-\frac{\lambda-\mu}{\lambda}\left[j(\delta\lambda)-\left(\begin{array}{c}j\\
2\end{array}\right)(\delta\lambda)^2\right.\\
&\ \ \ \left.+\ldots
+(-1)^i\left(\begin{array}{c}j\\
i\end{array}\right)(\delta\lambda)^{i+1} +\ldots
+(-1)^{j}(\delta\lambda)^{j}\right]\\
&{\buildrel{a.s}\over \leq}\overline{q}(t)
{\buildrel{a.s}\over \leq}1-\frac{\lambda-\mu}{\lambda}\left[(j-1)(\delta\lambda)-\left(\begin{array}{c}j-1\\
2\end{array}\right)(\delta\lambda)^2+\ldots\right.\\
&\ \ \ \left.+(-1)^i\left(\begin{array}{c}j-1\\
i\end{array}\right)(\delta\lambda)^{i+1} +\ldots
+(-1)^{j-1}(\delta\lambda)^{j-1}\right]
\end{aligned}
\end{equation}
for all $t\in[(j-1)\delta,j\delta)$.

Assume now that $\delta$ vanishes. Then taking
$j=\lfloor\frac{t}{\delta}\rfloor$, where $\lfloor \cdot\rfloor$
denotes the integer part of the number, and passing from the
sequence of sums to the limiting expression, for $\overline{q}(t)$
we arrive at
\begin{eqnarray}\label{1.11}
\overline{q}(t)&=&1-\frac{\lambda-\mu}{\lambda}\left[1-\lim_{\delta\downarrow0}
(1-\lambda\delta)^{\frac{t}{\delta}}\right]\nonumber\\
&=&1-\frac{\lambda-\mu}{\lambda}\left(1-\mathrm{e}^{-\lambda
t}\right),
\end{eqnarray}
for all $t\in[0, \infty)$.

Summarizing all of this we have as follows. Consider the interval
$[0,\infty)$, and build the system of partitions of this interval
$\{\Pi_\delta\}$. Then for any point $t\in(0,\infty)$ there is the
integer value $j_\delta$ (depending on $\delta$) such that $\delta
j_\delta\leq t<\delta j_\delta+\delta$. As $\delta$ vanishes and,
correspondingly, $j_\delta$ increases to infinity, then in the point
$t$ we have the convergence $\overline{q}(t){\buildrel
\mathrm{a.s.}\over=}\lim_{N\to\infty}\overline{q}_N(t)$, where
$\overline{q}(t)$ satisfies \eqref{1.11}. According to the
construction, this a.s. convergence is uniform on any bounded subset
of the interval $[0,\infty)$.

 Hence, in view of $q(t)=1-\overline{q}(t)$, \eqref{1.11} implies \eqref{1.10}, and
the desired relationship is proved.

\section{A general closed network containing several satellite
stations, one of which is bottleneck}\label{General}

Consider now a more general case of the model, where there are $k$
satellite stations numbered 1,2,\ldots,$k$, and only the $k$th
satellite station is bottleneck, while the first $k-1$ satellite
stations are non-bottleneck. In this case as in \cite{Abramov 2001}
we use the level crossing method. Note, that the approach of
\cite{Abramov 2001} based on straightforward generalization of the
level crossing method related to Markovian queues was mistaken. The
error of \cite{Abramov 2001} was corrected in another paper
\cite{Abramov 2006} and then in book \cite{Abramov 2009}, and here
we follow by the improved method. Let us recall this method in the
case of the standard $GI/M/1$ queueing system and then explain how
the method can be adapted to this network with bottleneck.

Consider the $GI/M/1$ queueing system, where an interarrival time
has the probability distribution function $G(x)$ with the
expectation $\frac{1}{\lambda}$, a service time is exponentially
distributed with parameter $\mu$, and $\lambda<\mu$. For a busy
period of this system, let $f(i)$ denote the number of cases during
that busy period where a customer at the moment of his/her arrival
finds exactly $i$ other customers in the system. Clearly, that
$f(0)=1$ with probability 1. Let $t_{i,1}$, $t_{i,2}$, \ldots,
$t_{i,f(i)}$ denote the time moments of these arrivals, and let
$s_{i,1}$, $s_{i,2}$,\ldots, $s_{i,f(i)}$ denote the service
completions during the aforementioned busy period, at which there
remain exactly $i$ customers in the system.

Consider the intervals
\begin{equation}\label{2.1}
[t_{i,1}, s_{i,1}), [t_{i,2}, s_{i,2}), \ldots, [t_{i,f(i)},
s_{i,f(i)})
\end{equation}
and
\begin{eqnarray}\label{2.2}
&&[t_{i+1,1}, s_{i+1,1}), [t_{i+1,2}, s_{i+1,2}), \ldots,\nonumber\\
&&[t_{i+1,f(i+1)}, s_{i+1,f(i+1)})
\end{eqnarray}
The intervals of (\ref{2.2}) all are contained in the intervals of
(\ref{2.1}). Let us delete the intervals of (\ref{2.2}) from those
of (\ref{2.1}) and merge the ends. Then we obtain a special type of
branching process $\{f(i)\}$, which according to its construction is
\textit{not} a standard branching process, because the number of
offspring generated by particles of different generations are not
independent random variables. (Only they are independent in the case
of the $M/M/1$ system.) However, the process satisfies the property:
$\mathrm{E}f(i)=\varphi^{i}$, $i=0,1,\ldots$, which is similar to
the respective property of a standard Galton-Watson branching
process. The branching process $\{f(i)\}$ is called $GI/M/1$
\textit{type branching process}. The properties of this branching
process are discussed in \cite{Abramov 2006} and \cite{Abramov
2009}.

The value $\varphi:=\mathrm{E}f(1)$ is calculated as follows. We
have the equation
$$
\mathrm{E}f(1)=\varphi=\sum_{i=0}^{\infty}\varphi^{i}\int_0^\infty\mathrm{e}^{-\mu
x}\frac{(\mu x)^{i}}{i!}\mathrm{d}G(x),
$$
where $\varphi$ is the least positive root of the functional
equation $z=\widehat{G}(\mu-\mu z)$; $\widehat{G}(s)$, $s\geq0$, is
the Laplace-Stieljies transform of the probability distribution
function $G(x)$.

Let us now consider the closed queueing network with $k$ satellite
stations. The first $k-1$ satellite stations are non-bottleneck,
i.e. the relation $\mu_j>\lambda p_j$ is satisfied, while the $k$th
satellite station is a bottleneck station, and then the relation
$\mu_k<\lambda p_k$ is satisfied (Recall that by $p_j$,
$j=1,2,\ldots,k$, we denote the routing probabilities.) The input
rate to the $j$th satellite station at time $t$ is $\lambda p_j
\mathcal{N}_t$, where the random variable $\mathcal{N}_t$
($\mathcal{N}_t\leq N$ with probability 1) is the number of units it
the hub at time $t$.

By the level crossing method one can study the queue-length process
in any satellite station $j$, $j=1,2,\ldots,k-1$. Consider first the
time interval [0, $\delta$), where $\delta$ is a fixed sufficiently
small value. Then, similarly to the construction in Section
\ref{OneSatellite} in the case of only one satellite station, for
$N$ increasing to infinity with probability approaching 1 the
queue-length in the hub at time moment $t=\delta$ is asymptotically
equal to $N[1-(\lambda p_k-\mu_k)\delta]$. If there are several
non-bottleneck satellite stations, then the queue-length in the hub
at time moment $\delta$ is asymptotically equal to the same
aforementioned value $N[1-(\lambda p_k-\mu_k)\delta]$ because the
queue-lengths in non-bottleneck satellite stations all are finite
with probability 1, and their contribution is therefore negligible.
More detailed arguments are as follows.


In any satellite station $j<k$ the length of a busy period is finite
with probability 1. Hence, as $N$ increases unboundedly, the number
of busy periods in the interval [0, $\delta$) increases to infinity
as well, and at time moment $\delta$ the total number of units in
all of satellite stations $j<k$ with probability 1 is bounded, i.e.
it is negligible compared to $N$. Therefore, the number of customers
in the hub is asymptotically evaluated as $N[1-(\lambda
p_k-\mu_k)\delta]$ as well.

For the $j$th satellite station, $j<k$, let us consider the last
busy period that finished before time moment $\delta$. Let
$f_j(i)$ denotes the number of cases during that busy period that
a customer at the moment of his/her arrival find $i$ other
customers in the system, let $t_{j,i,1}$,
$t_{j,i,2}$,\ldots,$t_{j,i,f_j(i)}$ be the moments of these
arrivals, and let $s_{j,i,1}$, $s_{j,i,2}$,\ldots,$s_{j,i,f_j(i)}$
be the moments of service completions that there remain exactly
$i$ units in the satellite station. We have the intervals
\begin{equation}\label{2.4}
[t_{j,i,1}, s_{j,i,1}), [t_{j,i,2}, s_{j,i,2}), \ldots,
[t_{j,i,f(i)}, s_{j,i,f_j(i)})
\end{equation}
and
\begin{eqnarray}\label{2.5}
&&[t_{j,i+1,1}, s_{j,i+1,1}), [t_{j,i+1,2}, s_{j,i+1,2}),
\ldots,\nonumber\\
&&[t_{j,i+1,f(i+1)}, s_{j,i+1,f_j(i+1)})
\end{eqnarray}
which are similar to the intervals (\ref{2.1}) and (\ref{2.2})
considered before. (The only difference in the additional index $j$
indicating the $j$th satellite station.) Delete the intervals
(\ref{2.5}) from those of (\ref{2.4}) and merge the ends. As $N$
increases unboundedly, the process $\{f_j(i)\}_{i\geq0}$ in a random
sequence in $i$. Each of the random variables $f_j(i)$ converges
a.s. (as $N\to\infty$) to the limiting random variable, which is the
number of offspring in the $i$th generation of the $GI/M/1$ type
branching process. So, we have the a.s. convergence to the $GI/M/1$
type branching process in the sense that for all $i$ this a.s.
convergence holds. Let us find $\mathrm{E}f_j(1)$. Similarly to the
above case of the $GI/M/1$ queueing system for sufficiently large
$N$ we have the equation
\begin{equation}\label{2.6}
\begin{aligned}
z&=\frac{p_j\widehat{F}_{\lfloor N[1-(\lambda
p_k-\mu_k)\delta]\rfloor}(\mu_j N-\mu_j
Nz)}{1-(1-p_j)\widehat{F}_{\lfloor N[1-(\lambda
p_k-\mu_k)\delta]\rfloor}(\mu_j N-\mu_j Nz)}[1+o(1)]\\
&=\frac{p_j\widehat{G}_{\lfloor N[1-(\lambda
p_k-\mu_k)\delta]\rfloor}\left(\frac{\mu_j N -\mu_j Nz}{\lfloor
N[1-(\lambda
p_k-\mu_k)\delta]\rfloor}\right)}{1-(1-p_j)\widehat{G}_{\lfloor
N[1-(\lambda p_k-\mu_k)\delta]\rfloor}\left(\frac{\mu_j N-\mu_j
Nz}{\lfloor N[1-(\lambda p_k-\mu_k)\delta]\rfloor}\right)}[1+o(1)],
\end{aligned}
\end{equation}
where $\widehat{F}_N(s)$ denotes the Laplace-Stieltjes transform of
the probability distribution function $F_N(x)=G_N(Nx)=G(Nx)$. For
details of the derivation of \eqref{2.6} see the Appendix. As
$N\to\infty$, in limit we have the equation
\begin{equation}\label{2.7}
z=\frac{p_j\widehat{G}\left(\frac{\mu_j -\mu_j z}{1-(\lambda
p_k-\mu_k)\delta}\right)}{1-(1-p_j)\widehat{G}\left(\frac{\mu_j
-\mu_j z}{1-(\lambda p_k-\mu_k)\delta}\right)}.
\end{equation}
Considering now the interval [$\delta$, $2\delta$), in the endpoint
$2\delta$, due to the arguments of Section \ref{OneSatellite} and
the arguments above in this section, the queue-length is
asymptotically equal to $N\frac{\lambda p_k-\mu_k}{\lambda
p_k}[2(\lambda p_k\delta)-(\lambda p_k\delta)^{2}]$. Therefore,
similarly to (\ref{2.7}) for the $j$th satellite station, $j<k$, we
have the equation
\begin{equation}\label{2.8}
z=\frac{p_j\widehat{G}\left(\frac{\mu_j -\mu_j z}{1-\frac{\lambda
p_k-\mu_k}{\lambda p_k}[2\lambda p_k\delta- (\lambda
p_k\delta)^{2}]}\right)}{1-(1-p_j)\widehat{G}\left(\frac{\mu_j
-\mu_j z}{1-\frac{\lambda p_k-\mu_k}{\lambda p_k}[2\lambda
p_k\delta- (\lambda p_k\delta)^{2}]}\right)}.
\end{equation}
In an arbitrary interval [$(i-1)\delta$, $i\delta$), in its endpoint
$i\delta$ we correspondingly have the equation:
\begin{equation}\label{2.9}
z=\frac{p_j\widehat{G}\left(\frac{\mu_j -\mu_j
z}{U(i,\delta)}\right)}{1-(1-p_j)\widehat{G}\left(\frac{\mu_j -\mu_j
z}{U(i,\delta)}\right)},
\end{equation}
where
\begin{eqnarray}
U(i,\delta)&=&
1-\frac{\lambda p_k-\mu_k}{\lambda p_k}[i\delta\lambda p_k-\left(\begin{array}{c}i\\
2\end{array}\right)
(\lambda p_k\delta)^{2}+\ldots\nonumber\\ &&+(-1)^{l}\left(\begin{array}{c}i\\
l\end{array}\right)(\lambda p_k\delta)^{l}+\ldots\nonumber\\
&&+(-1)^{i}(\lambda p_k\delta)^i].\nonumber
\end{eqnarray}
Taking $i=\lfloor\frac{t}{\delta}\rfloor$, in the limit as
$\delta\to0$ we obtain:
\begin{equation}\label{2.10}
z=\frac{p_j\widehat{G}\left(\frac{\mu_j -\mu_j
z}{\overline{q}(t)}\right)}{1-(1-p_j)\widehat{G}\left(\frac{\mu_j
-\mu_j z}{\overline{q}(t)}\right)},
\end{equation}
where $\overline{q}(t)=1-\frac{\lambda p_k-\mu_k}{\lambda
p_k}(1-\mathrm{e}^{-\lambda p_kt})$.

Thus, denoting by $\varphi_j(t)$ the root of the functional equation
(\ref{2.10}), we arrive at the following statement.

\medskip

\begin{thm}\label{thm0} The queue-length distribution in the
non-bottleneck satellite station $j<k$ is
\begin{eqnarray}\label{2.11}
&&\mathrm{P}\{Q_j(t)=i\}=\rho_j(t)[\varphi_j(t)]^{i-1}[1-\varphi_j(t)],\nonumber\\
&&i=1,2,\ldots,
\end{eqnarray}
where
$$
\rho_j(t)= \frac{\lambda \overline{q}(t) p_j}{\mu_j},
$$
$$
\overline{q}(t)=1-\frac{\lambda p_k-\mu_k}{\lambda
p_k}(1-\mathrm{e}^{-\lambda p_kt}),
$$
and $\varphi_j(t)$ the root of the functional equation (\ref{2.10}).
\end{thm}

\section{Continuity of the queue-length processes}\label{Continuity}

\subsection{Formulation of the main results}\label{Main}
\begin{thm}\label{thm1} For small positive $\epsilon$ assume that Condition
\eqref{0.2} is satisfied. Then, for $\varphi_j(t)$ in \eqref{2.11}
the following bounds are true:
$$
\rho_j(t)-\epsilon_1(t)\leq\varphi_j(t)\leq\rho_j(t)+\epsilon_2(t),
$$
where
$$
\epsilon_1(t)=\min\left\{\rho_j(t)-\ell_j(t), \
\frac{2\epsilon}{p_j}[1-\ell_j(t)]\right\},
$$
$$
\epsilon_2(t)=\min\left\{1+\frac{[a_j(t)]^{2}}{b_j(t)}[\ell_j(t)-1]-\rho_j(t),
\ \frac{2\epsilon}{p_j}[1-\ell_j(t)]\right\},
$$
$\ell_j(t)$ is the least root of the equation
$$
z=\mathrm{e}^{-a_j(t)\mu_j+a_j(t)\mu_j z},
$$
\begin{equation*}
a_j(t)=\frac{1}{\lambda p_j\overline{q}(t)},
\end{equation*}
and
\begin{equation*}
\begin{aligned}
b_j(t)&=\frac{1}{p_j[\overline{q}(t)]^2}\left(r-\frac{1}{\lambda^{2}}\right)
+\left(\frac{1}{1-p_j}+\frac{1-2p_j}{p_j^{2}}\right)\frac{1}{\lambda^{2}[\overline{q}(t)]^{2}}\\
&\ \ \ +\frac{1}{\lambda^{2}p_j^{2}[\overline{q}(t)]^{2}}.
\end{aligned}
\end{equation*}
\end{thm}

\begin{thm}\label{thm2} Under the assumptions made in Theorem \ref{thm1}
assume additionally that the probability distribution function
$G(x)$ belongs either to the class NBU or to the class NWU. Then,
for $\varphi_j(t)$ in \eqref{2.11} the following bounds are true:
$$
\rho_j(t)-\epsilon_3(t)\leq\varphi_j(t)\leq\rho_j(t)+\epsilon_4(t),
$$
where
$$
\epsilon_3(t)=\min\left\{\rho_j(t)-\ell_j(t), \
\frac{\epsilon}{p_j}[1-\ell_j(t)]\right\},
$$
$$
\epsilon_4(t)=\min\left\{1+\frac{[a_j(t)]^{2}}{b_j(t)}[\ell_j(t)-1]-\rho_j(t),
\ \frac{\epsilon}{p_j}[1-\ell_j(t)]\right\},
$$
and $\ell_j(t)$, $a_j(t)$ and $b_j(t)$ are as in Theorem \ref{thm1}.
\end{thm}

Note, that the difference between $\epsilon_1(t)$ (or
$\epsilon_2(t)$) \ in Theorem \ref{thm1} and $\epsilon_3(t)$ (or
$\epsilon_4(t)$) in Theorem \ref{thm2} is only in one term in the
$\min$ function. While $\epsilon_1(t)$ (or $\epsilon_2(t)$) contains
the term $\frac{2\epsilon}{p_j}[1-\ell_j(t)]$, the corresponding
term in $\epsilon_3(t)$ (or $\epsilon_4(t)$) is
$\frac{\epsilon}{p_j}[1-\ell_j(t)]$.

\subsection{Background derivations and the proof of the theorems}
The results of continuity that used to prove Theorems \ref{thm1} and
\ref{thm2} are based on a recent result obtained in \cite{Abramov
cont}. Recall it. Let $G_1(x)$ and $G_2(x)$ be two probability
distribution functions belonging to the class $\mathcal{G}(a,b)$ of
probability distributions functions of positive random variables
having the first and second moments $a$ and $b$ respectively,
$b>a^2$. Assume that
\begin{equation}\label{e1}
\sup_{x>o}|G_1(x)-G_2(x)|<\kappa,
\end{equation}
where $\kappa<1-\frac{{a}^2}{b}$. Let $\widehat{G}_1(s)$ and
$\widehat{G}_2(s)$ be the corresponding Laplace-Stieltjes transforms
of $G_1(x)$ and $G_2(x)$ and let $\gamma_{G_1}$ and $\gamma_{G_2}$
be the least positive roots of corresponding functional equations
$z=\widehat{G}_1(\mu-\mu z)$ and $z=\widehat{G}_2(\mu-\mu z)$, where
$\mu>\frac{1}{a}$ is some real number. Then
$|\gamma_{G_1}-\gamma_{G_2}|<\kappa(1-\ell)$, where $\ell$ is the
least root of the equation
$$
x=\mathrm{e}^{-a\mu+a\mu x}.
$$
On the other hand, according to the results of Rolski \cite{Rolski
1972}, the guaranteed bounds for any of probability distribution
functions $G_1(x)$ and $G_2(x)$ (i.e. not necessarily satisfying
\eqref{e1}) having the first two moments $a$ and $b$ respectively
are given by
\begin{equation}\label{b1}
\ell\leq\gamma_{G_i}\leq 1+\frac{a^{2}}{b}(\ell-1), \ i=1,2.
\end{equation}

In our case, \eqref{2.10} can be formally rewritten as
$$
z=\widehat{H}_j(\mu_j-\mu_j z),
$$
where
$$
\widehat{H}_j(s)=\frac{p_j\widehat{G}\left(\frac{s}{\overline{q}(t)}\right)}
{1-(1-p_j)\widehat{G}\left(\frac{s}{\overline{q}(t)}\right)}
$$
is the Laplace-Stieltjes transform of a positive random variable
having the probability distribution $H_j(x)$,
$$
H_j(x)=\sum_{i=1}^\infty p_j(1-p_j)^{i-1}G^{*i}(\overline{q}(t)x),
$$
$G^{*i}(x)$ denotes $i$-fold convolution of the probability
distribution function $G(x)$ with itself, and
$\overline{q}(t)=1-\frac{\lambda p_k-\mu_k}{\lambda
p_k}(1-\mathrm{e}^{-\lambda p_kt})$. Our challenge now is to find
the bounds for the least positive root $\gamma_{H_j}$ of the
functional equation in \eqref{2.10} similar to those given in
\eqref{b1} and then, applying the aforementioned result of
\cite{Abramov cont}, find the continuity bounds for non-stationary
distributions in non-bottleneck stations. (The notation
$\gamma_{H_j}$ is used here (along with the other notation
$\varphi_j(t)$ in the formulations of Theorems \ref{thm0},
\ref{thm1} and \ref{thm2}) for consistency with the notation such as
$\gamma_{G_1}$ and $\gamma_{G_2}$ that introduced before.)

The probability distribution function $G(x)$ is assumed to have the
expectation $\frac{1}{\lambda}$ and the second moment
$r>\frac{1}{\lambda^{2}}$. Let us find the expectation and second
moment of the probability distribution $H_j(x)$. The best way for
deriving these numerical characteristics is to use Wald's identities
as follows. Let $\xi_1$, $\xi_2$, \ldots be a sequence of
independent and identically distributed random variables, let $\tau$
be an integer random variable independent of the sequence $\xi_1$,
$\xi_2$, \ldots. Denote $S_\tau=\xi_1+\xi_2+\ldots+\xi_\tau$. Then,
\begin{eqnarray}
\mathrm{E}S_\tau&=&\mathrm{E}\xi_1\mathrm{E}\tau,\label{Wald0}\\
\mathrm{Var}(S_\tau)&=&\mathrm{Var}(\xi_1)\mathrm{E}\tau+\mathrm{Var}(\tau)(\mathrm{E}\xi_1)^2.\label{Wald}
\end{eqnarray}
In our case, the random variable $\tau$ is a geometrically
distributed random variable,
$\mathrm{Pr}\{\tau=n\}=p_j(1-p_j)^{n-1}$, $n\geq1$. Therefore,
$$\mathrm{E}\tau=\frac{1}{p_j},$$ and $$\mathrm{Var}(\tau)=\frac{1}{1-p_j}+\frac{1-2p_j}{p_j^{2}}.$$
The random variable $\xi_1$ has the probability distribution
function $G(\overline{q}(t)x)$. Hence,
$$\mathrm{E}\xi_1=\frac{1}{\lambda \overline{q}(t)},$$ and
$$\mathrm{Var}(\xi_1)=\frac{1}{[\overline{q}(t)]^2}\left(r-\frac{1}{\lambda^{2}}\right).$$

So, according to Wald's equations \eqref{Wald0} and \eqref{Wald} we
obtain:
\begin{eqnarray}
\int_0^\infty x\mathrm{d}H_j(x)&=&\frac{1}{\lambda
p_j\overline{q}(t)},\label{b0}\\
\int_0^\infty
x^2\mathrm{d}H_j(x)&=&\frac{1}{p_j[\overline{q}(t)]^2}\left(r-\frac{1}{\lambda^{2}}\right)
+\left(\frac{1}{1-p_j}+\frac{1-2p_j}{p_j^{2}}\right)\frac{1}{\lambda^{2}[\overline{q}(t)]^{2}}\label{b2}\\
&&+\frac{1}{\lambda^{2}p_j^{2}[\overline{q}(t)]^{2}}.\nonumber
\end{eqnarray}

\textit{Proof of Theorem \ref{thm1}.} It follows from \eqref{b0} and
\eqref{b2} and the aforementioned result by Rolski \cite{Rolski
1972} that $\gamma_{H_j}$ satisfies the inequalities
\begin{equation}\label{d1}
\ell_j(t)\leq\gamma_{H_j}\leq1+\frac{[a_j(t)]^{2}}{b_j(t)}[\ell_j(t)-1],
\end{equation}
where $\ell_j(t)$ is the least root of the equation
\begin{equation}\label{d2}
z=\mathrm{e}^{-a_j(t)\mu_j+a_j(t)\mu_j z},
\end{equation}
\begin{equation}\label{d3}
a_j(t)=\frac{1}{\lambda p_j\overline{q}(t)},
\end{equation}
and
\begin{equation}\label{d4}
\begin{aligned}
b_j(t)&=\frac{1}{p_j[\overline{q}(t)]^2}\left(r-\frac{1}{\lambda^{2}}\right)
+\left(\frac{1}{1-p_j}+\frac{1-2p_j}{p_j^{2}}\right)\frac{1}{\lambda^{2}[\overline{q}(t)]^{2}}\\
&\ \ \ +\frac{1}{\lambda^{2}p_j^{2}[\overline{q}(t)]^{2}}.
\end{aligned}
\end{equation}

Assume that \eqref{0.2} is satisfied. Then according to the
characterization theorem of Azlarov and Volodin \cite{Azlarov and
Volodin 1987}, \cite{Abramov 2008b} we have \eqref{0.5}, and then
\begin{equation}\label{c1}
\begin{aligned}
\sup_{s\geq0}\left|\widehat{G}\left(\frac{s}{\overline{q}(t)}\right)-\frac{\lambda
\overline{q}(t)}{\lambda \overline{q}(t)+s}\right|&=
\sup_{s>0}\left|\int_0^\infty\mathrm{e}^{-sx}\mathrm{d}G(\overline{q}(t)x)-\int_0^\infty
\mathrm{e}^{-sx}\mathrm{d}(1-\mathrm{e}^{-\lambda \overline{q}(t) x})\right|\\
&\leq\sup_{s>0}\int_0^\infty
s\mathrm{e}^{-sx}\sup_{x\geq0}|G(\overline{q}(t)x)-(1-\mathrm{e}^{-\lambda
\overline{q}(t)x})|\mathrm{d}x\\
&=\sup_{s>0}\int_0^\infty
s\mathrm{e}^{-sx}\underbrace{\sup_{x\geq0}|G(x)-(1-\mathrm{e}^{-\lambda
x})|}_{\leq 2\epsilon}\mathrm{d}x\\
&\leq2\epsilon.
\end{aligned}
\end{equation}
In turn, from \eqref{c1} we obtain
\begin{equation}\label{c2}
\begin{aligned}
&\sup_{s\geq0}\left|\widehat{H}_j(s)-\frac{\lambda
p_j\overline{q}(t)}{\lambda p_j\overline{q}(t)+s}\right|\\&=
\sup_{s\geq0}\left|\sum_{i=1}^\infty
p_j(1-p_j)^{i-1}\left[\widehat{G}^i\left(\frac{s}{\overline{q}(t)}\right)-\left(\frac{\lambda
\overline{q}(t)}{\lambda \overline{q}(t)+s}\right)^i\right]\right|\\
&\leq\sum_{i=1}^\infty
p_j(1-p_j)^{i-1}\underbrace{\sup_{s\geq0}\left|\widehat{G}^i\left(\frac{s}{\overline{q}(t)}\right)-\left(\frac{\lambda
\overline{q}(t)}{\lambda \overline{q}(t)+s}\right)^i\right|}_{\leq
2i\epsilon}\\
&\leq2\epsilon\sum_{i=1}^\infty
ip_j(1-p_j)^{i-1}\\
&=\frac{2\epsilon}{p_j}.
\end{aligned}
\end{equation}

Then, the results of \cite{Abramov cont} enables us to conclude as
follows. Let us consider the functional equation
\begin{equation}\label{e5}
z=\frac{\lambda p_j\overline{q}(t)}{\lambda
p_j\overline{q}(t)+\mu_j-\mu_jz}.
\end{equation}
The right-hand side of \eqref{e5} can be written as
$\widehat{\Pi}(\mu_j-\mu_j z)$, where $\Pi(s)$ is the
Laplace-Stieltjes transform of exponential distribution with the
parameter $\lambda p_j\overline{q}(t)$

The least positive root of equation \eqref{e5} is
$$
\rho_j(t)=\frac{\lambda p_j\overline{q}(t)}{\mu_j}.
$$
Apparently,
\begin{equation}\label{d5}
\ell_j(t)\leq\gamma_{H_j}\leq1+\frac{[a_j(t)]^{2}}{b_j(t)}[\ell_j(t)-1],
\end{equation}
where $\ell_j(t)$, $a_j(t)$ and $b_j(t)$ are defined in
\eqref{d2}-\eqref{d4}. The lower and upper bounds in \eqref{d5} are
natural bounds obtained from the results by Rolski \cite{Rolski
1972}.

Along with \eqref{d5} we have also the following inequality
\begin{equation}\label{f2}
\ell_j(t)\leq\rho_{j}(t)\leq1+\frac{[a_j(t)]^{2}}{b_j(t)}[\ell_j(t)-1],
\end{equation}
which holds true because of Condition \eqref{f1}.

On the other hand, following the results in \cite{Abramov cont} and
\eqref{c2},
$$
|\gamma_{H_j}-\rho_j(t)|\leq\frac{2\epsilon}{p_j}[1-\ell_j(t)].
$$
So,
\begin{equation}\label{d6}
\rho_j(t)-\frac{2\epsilon}{p_j}[1-\ell_j(t)]\leq\gamma_{H_j}\leq\rho_j(t)-\frac{2\epsilon}{p_j}[1-\ell_j(t)].
\end{equation}
Amalgamating \eqref{d5} and \eqref{d6} we arrive at
the following bounds:
\begin{equation}\label{d7}
\rho_j(t)-\epsilon_1(t)\leq\gamma_{H_j}\leq\rho_j(t)+\epsilon_2(t),
\end{equation}
where
$$
\epsilon_1(t)=\min\left\{\rho_j(t)-\ell_j(t), \
\frac{2\epsilon}{p_j}[1-\ell_j(t)]\right\},
$$
and
$$
\epsilon_2(t)=\min\left\{1+\frac{[a_j(t)]^{2}}{b_j(t)}[\ell_j(t)-1]-\rho_j(t),
\ \frac{2\epsilon}{p_j}[1-\ell_j(t)]\right\}.
$$
The theorem is proved.

\medskip
\textit{Proof of Theorem \ref{thm2}.} The proof of this theorem is
similar to that of Theorem \ref{thm1}. The only difference that
instead of \eqref{0.5} given by characterization theorem of Azlarov
and Volodin \cite{Azlarov and Volodin 1987}, \cite{Abramov 2008b}
one should use \eqref{0.7} from \cite{Abramov 2008b}. In this case,
instead of estimate \eqref{c1} we arrive at
\begin{equation*}\label{c1*}
\begin{aligned}
\sup_{s\geq0}\left|\widehat{G}\left(\frac{s}{\overline{q}(t)}\right)-\frac{\lambda
\overline{q}(t)}{\lambda \overline{q}(t)+s}\right|\leq\epsilon,
\end{aligned}
\end{equation*}
and then instead of \eqref{c2} at
\begin{equation*}\label{c2*}
\begin{aligned}
\sup_{s\geq0}\left|\widehat{H}_j(s)-\frac{\lambda
p_j\overline{q}(t)}{\lambda
p_j\overline{q}(t)+s}\right|\leq\frac{\epsilon}{p_j}.
\end{aligned}
\end{equation*}
Then instead of \eqref{d6} we have
\begin{equation*}\label{d6*}
\rho_j(t)-\frac{\epsilon}{p_j}[1-\ell_j(t)]\leq\gamma_{H_j}\leq\rho_j(t)-\frac{\epsilon}{p_j}[1-\ell_j(t)],
\end{equation*}
and instead of \eqref{d7} we obtain
\begin{equation*}\label{d7*}
\rho_j(t)-\epsilon_3(t)\leq\gamma_{H_j}\leq\rho_j(t)+\epsilon_4(t),
\end{equation*}
with
$$
\epsilon_3(t)=\min\left\{\rho_j(t)-\ell_j(t), \
\frac{\epsilon}{p_j}[1-\ell_j(t)]\right\},
$$
and
$$
\epsilon_4(t)=\min\left\{1+\frac{[a_j(t)]^{2}}{b_j(t)}[\ell_j(t)-1]-\rho_j(t),
\ \frac{\epsilon}{p_j}[1-\ell_j(t)]\right\}.
$$
The theorem is proved.
\bigskip

\appendix{\sc Appendix: Derivation of \eqref{2.6}}

Let $\tau_1$, $\tau_2$, \ldots be a sequence of service completions
in the hub, and let $Q_0(\tau_1)$, $Q_0(\tau_2)$, \ldots the
queue-lengths in the hub in these time instants $\tau_1$, $\tau_2$,
\ldots. Then, the probability distribution function of the first
interarrival time in the satellite node $j$ can be represented as
$$
H_{j,1}(x)=\sum_{i=0}^\infty
p_j(1-p_j)^iF_{r_0}*F_{r_1}*\ldots*F_{r_i}(x),\leqno{(A.1)}
$$
where $F_{r_0}(x)=F_{N}(x)=G_{N}(Nx)$, $F_{r_l}(x)=\sum_{u=1}^\infty
\mathrm{Pr}\{Q_0(\tau_l)=u\}G_u(ux)$, $l\geq1$, the asterisk denotes
convolution between the probability distribution functions.

Denote the sequence of interarrival time distributions in the
satellite node $j$ by $H_{j,1}(x)$, $H_{j,2}(x)$,\ldots. Let us
derive the representation for $H_{j,2}(x)$. Keeping in mind that the
time instants $\tau_1$, $\tau_2$,\ldots, $\tau_i$, \ldots can be the
moments of the first arrival to the satellite node $j$ with
probabilities $p_j$, $p_j(1-p_j)$,\ldots, $p_j(1-p_j)^{i-1}$,\ldots
respectively, we obtain:
$$
H_{j,2}(x)=\sum_{l=1}^\infty\sum_{i=0}^\infty p_j^{2}(1-p_j)^{i+l-1}
F_{r_l}*F_{r_{l+1}}*\ldots*F_{r_{l+i}}(x),\leqno{(A.2)}
$$
and recurrently,
$$
\begin{aligned}
H_{j,v}(x)&=\sum_{k_1=1}^\infty\cdots\sum_{k_{v-1}=1}^\infty\sum_{i=0}^\infty
p_j^v(1-p_j)^{i+k_1+k_2+\ldots+k_{v-1}-v+1}\\
&\times\left[F_{r_{k_1+\ldots+k_{v-1}}}*F_{r_{k_1+\ldots+k_{v-1}+1}}*\ldots*
F_{r_{k_1+\ldots+k_{v-1}+i}}(x)\right].
\end{aligned}\leqno{(A.3)}
$$
Let us now use the assumption given by \eqref{0.1add}. According to
this assumption, we have $F_{i}(x)\leq F_N(x)$ for all $i\leq N$.
Hence, for the convolution given in the right-hand side of (A.3) we
have:
$$
\left[F_{r_{k_1+\ldots+k_{v-1}}}*F_{r_{k_1+\ldots+k_{v-1}+1}}*\ldots*
F_{r_{k_1+\ldots+k_{v-1}+i}}(x)\right]\leq F_N^{*i+1}(x),
$$
where $F_N^{*i+1}(x)$ is $(i+1)$-fold convolution of the probability
distribution function $F_N(x)$ with itself. Therefore,
$$
H_{j,v}(x)\leq\sum_{i=0}^\infty
p_j(1-p_j)^iF_N^{*i+1}(x)\leqno{(A.4)}
$$
for all $v\geq1$ and $x\geq0$, and moreover formula (A.4) remains
correct when $v$ is a positive integer random variable not smaller
than 1, due to the formula for the total probability.

Let us consider the intervals of the type of \eqref{2.4} and
\eqref{2.5} for the node $j$. To this end, consider the procedure of
deleting the intervals of \eqref{2.5} from those \eqref{2.4} and
merging the ends. Denote also $N_j^+$ and $N_j^-$ for maximum and
minimum numbers of customers in the hub during the busy period in
the $j$th satellite station where these intervals are considered.
Clearly, that $N_j^+\leq N$, and $N_j^-\geq N-\sum_{i=1}^k n_i$,
where $n_i$ denote the number of unit arrivals to the satellite node
$i$ from the hub. Let $\nu_j$ denote the number of units served
during the busy period in the satellite station $j$. The
distribution of $\nu_j$ generally depends on the state of the
network at the moment of the busy period start, and it is rightly to
write $\nu_j(\mathcal{S})$. This state $\mathcal{S}$ includes the
number of units in the hub, which is assumed to be large enough such
that as $N\to\infty$ it becomes asymptotically equivalent to $N$.

Apparently, for the expected number of units served during the busy
period we have the inequality $\mathrm{E}\nu_j(\mathcal{S})\leq C$,
where $C$ is the expected number of served units during a busy
period of the $GI/M/1$ queueing system with independently and
identically distributed interarrival times having the probability
distribution $\sum_{i=0}^\infty p_j(1-p_j)^iF_N^{*i+1}(x)$ and
exponentially distributed service time with parameter $\mu_jN$. (The
mentioned probability distribution of a service time is given by the
right-hand side of (A.4).) For a non-bottleneck station $j$ the
value $C$ is finite, since the expected interarrival time associated
with the probability distribution function $\sum_{i=0}^\infty
p_j(1-p_j)^iF_N^{*i+1}(x)$, which is equal to $\frac{1}{\lambda p_j
N}$, is assumed to be greater than $\frac{1}{\mu_j N}$. Hence, the
random sum $\sum_{i=1}^k n_i$ is a finite random variable with the
finite expectation. Therefore, as $N$ increases to infinity,
$\frac{N_j^-}{N}$ converges a.s. to 1. This means that, as
$N\to\infty$, the random variable $\nu_j$ converges in distribution
to the number of served units during a busy period of the $GI/M/1$
queueing system with independently and identically distributed
interarrival times having the probability distribution function
$\sum_{i=0}^\infty p_j(1-p_j)^iF_N^{*i+1}(x)$ and exponentially
distributed service times with parameter $\mu_jN$. As well, let
$f_{j}(\mathcal{S},n)$ denote the number of cases during the busy
period $\nu_j(\mathcal{S})$, when an arriving unit meets $n$ units
in the $j$th satellite station. As $N\to\infty$, the process
$\{f_{j}(\mathcal{S},n)\}$ converges in distribution to the
associated $GI/M/1$-type branching process. Specifically, let
$\varphi_{j,N}$ be the least positive root of the functional
equation
$$
\begin{aligned}
z_N&=\sum_{i=0}^\infty
p_j(1-p_j)^i\widehat{F}_N^{i+1}(\mu N-\mu Nz)\\
&=\frac{p_j\widehat{F}_N(\mu_jN-\mu_jNz)}{1-(1-p_j)\widehat{F}_N(\mu_jN-\mu_jNz)}.
\end{aligned}
$$
Then,
$$
\begin{aligned}
\lim_{N\to\infty}\mathrm{E}f_j(\mathcal{S},1)&=\lim_{N\to\infty}z_N\\
&=\lim_{N\to\infty}\frac{p_j\widehat{F}_N(\mu_jN-\mu_jNz)}{1-(1-p_j)\widehat{F}_N(\mu_jN-\mu_jNz)}.
\end{aligned}\leqno{(A.5)}
$$
Taking into account, that
$\widehat{F}_N(s)=\widehat{G}_N\left(\frac{s}{N}\right)$ and
$\lim_{N\to\infty}\widehat{G}_N(s)=G(s)$, from (A.5) we obtain
$$
\lim_{N\to\infty}\mathrm{E}f_j(\mathcal{S},1)=\varphi_j=
\frac{p_j\widehat{G}(\mu_j-\mu_j\varphi_j)}{1-(1-p_j)\widehat{G}(\mu_j\varphi_j-\mu_j\varphi_j)},
$$
where $\varphi_j=\lim_{N\to\infty}\varphi_{j,N}$.

Let us now derive \eqref{2.6}. Let $[0,\delta)$ be a small time
interval, and let $Z$ be the last busy period in this time interval
in the satellite station $j$. As $N$ becomes large enough, at the
end of the interval $[0,\delta)$ the number of units in the hub is
proportional to $[1-(\lambda p_k-\mu_k)\delta]N$. The following
arguments in the derivation of \eqref{2.6} are the same as above.
The only difference is that the probability distribution function
$F_N(x)$ is to be replaced with $F_{\lfloor [1-(\lambda
p_k-\mu_k)\delta]N\rfloor}(x)$, and, consequently, its
Laplace-Stieltjes transform $\widehat{F}_N(s)$ is to be replaced
with $\widehat{F}_{\lfloor [1-(\lambda
p_k-\mu_k)\delta]N\rfloor}(s)$. Therefore, instead of asymptotic
relation (A.5), we obtain the following asymptotic relation
$$
\begin{aligned}
\mathrm{E}f_j(\mathcal{S},1)&= \frac{p_j\widehat{F}_{\lfloor
N[1-(\lambda
p_k-\mu_k)\delta]\rfloor}(\mu_jN-\mu_jNz)}{1-(1-p_j)\widehat{F}_{\lfloor
N[1-(\lambda p_k-\mu_k)\delta]\rfloor}(\mu_jN-\mu_jNz)}[1+o(1)],
\end{aligned}
$$
which coincides with the asymptotic formula in \eqref{2.6}.

\end{document}